\documentclass[a4paper,12pt,oneside]{article}

\usepackage{fancyhdr}
\usepackage{epsfig}
\usepackage{amssymb,amsmath}
\usepackage{latexsym}
\usepackage{amsthm}

\newtheorem{theorem}{{\bf Theorem}}
\newtheorem{proposition}{{\bf Proposition}}

\newcommand{\R}{{\mathbb R}}
\newcommand{\pp}{{\bf p}}

\newcommand{\lr}{\leftrightarrow}

\begin{document}
{\bf Geometry of isophote curves.}
\vskip 3truemm
 Andr\'e DIATTA and Peter J. GIBLIN \footnote{ \footnotesize
 University of Liverpool. Department of Mathematical Sciences. M$\&$O Building, Peach Street, Liverpool, L69 7ZL, UK. adiatta@liv.ac.uk, pjgiblin@liv.ac.uk
\newline\noindent
$\bullet$ The first author was supported by the IST Programme of the European Union (IST-2001-35443).
\newline\noindent
$\bullet$ This work is a part of the DSSCV project supported by the IST Programme of the European Union (IST-2001-35443).
\\
The authors are very grateful to Prof. Mads Nielsen and Dr. Aleksander Pukhlikov for very fruitful discussions.
 }

\begin{abstract} We consider the intensity surface of a 2D image, we study the evolution of the symmetry sets (and medial axes) of
1-parameter families of iso-intensity curves.  This extends the investigation done on 1-parameter families of smooth plane curves
(Bruce and Giblin, Giblin and Kimia, etc.) to the general case when the family of curves includes a singular member, as will happen
if the curves are obtained by taking plane sections of a smooth surface, at the moment when the plane becomes tangent to the surface.
\end{abstract}

\section{Introduction}
Image data is often thought of as a collection of pixel values $I: Z^2\mapsto
Z_+$. The physical information is better captured by embedding the pixel values
in the real plane, as the pixelation and quantization are artifacts of the
camera, hence $I:{\bf R}^2\mapsto {\bf R}_+$.  The geometrical information of an image is
even better captured looking at the level sets $I(x)=I_0$, for all $I_0 \in {\bf R}_+$, that is, looking at the isophote curves  of the image.

Shape analysis using point-based representations or medial representations (such
as skeletons) has been widely applied on an object level demanding object
segmentation from the image data.
We propose to combine the object representation using a skeleton or symmetry set
representation and the appearance modelling by representing image information as
a collection of medial representations for the level-sets of an image.  As the level $I_0$  changes, the curves change like
sections of a smooth surface by parallel planes.

The qualitive changes in the medial representation of families of isophotes fall into two types: (1) those for which the isophotes
remain nonsingular  (see for example~\cite{growth-motion,IJCV03}) and (2) those for which one isophote at least is singular.
The symmetry set (SS) of a plane curve is the closure of the set of centres of circles which are tangent to the curve at two
or more different places.
The medial axis (MA) is the subset of the SS consisting of the closure of the locus of centres of circles which are maximal,
(maximal means that the minimum distance from the centre to the curve equals the radius).
Our aim is to extend the investigation to the case (2) when the family includes  singular curves, as is the case when one of
 the plane sections is tangent to the surface so that this section is a singular curve.
The final goal is to represent image smooth surfaces by the collection of all medial reprentations of isophotes, forming a
singular surface in scale space.

\smallskip

In this article, which is theoretical in nature, we work with the full SS, and  consider the
transitions which occur in the SS of a family
of plane sections of a  generic smooth surface in 3-space, as the plane moves through a position where it is tangent to the surface.
We investigate the local geometry of these families of curves and track the evolution of some crucial features of the SS and MA. In particular, we will trace and classify the patterns of  some special points, on the sections of a surface
 as the section passes through a tangential point, such as {\em vertices} (maxima and minima of curvature), {\em inflexions},
{\em triples of points  where a circle is tritangent} and the pattern of the centre of such a circle, {\em paires of points where
a circle is bitangent with a higher order contact at one of them}, etc. The vertices are crucial to the understanding of  the
 SS since it has branches which end at the centres of curvature at vertices. From the way in which vertices behave
 we can deduce a good deal about the evolution of the SS and its local number of branches. The {\em inflexions}
correspond to where the evolute of the curve, recedes to infinity. We also classify all possible scenarii of how vertices and
inflexions are distributed along the level curves.

Last, we produce examples of SS and MA illustrating the cases.

We are concerned with the local behaviour of symmetry sets (SS) and medial axis (MA) of plane sections
of generic\footnote{\footnotesize The genericity conditions will vary from case to case. See
\cite{diatta-giblin2004}.} smooth surfaces so we may assume that our surface
$M$ is given by an equation $z=f(x,y)$ for a smooth function $f$, which will often be assumed to be
 a polynomial of sufficiently high degree.
We shall take $M$ in Monge form, that is $f, f_x$ and $f_y$ all vanish at (0,0).

\begin{figure}
  \begin{center}
  \leavevmode
\epsfysize=1.2in
\epsffile{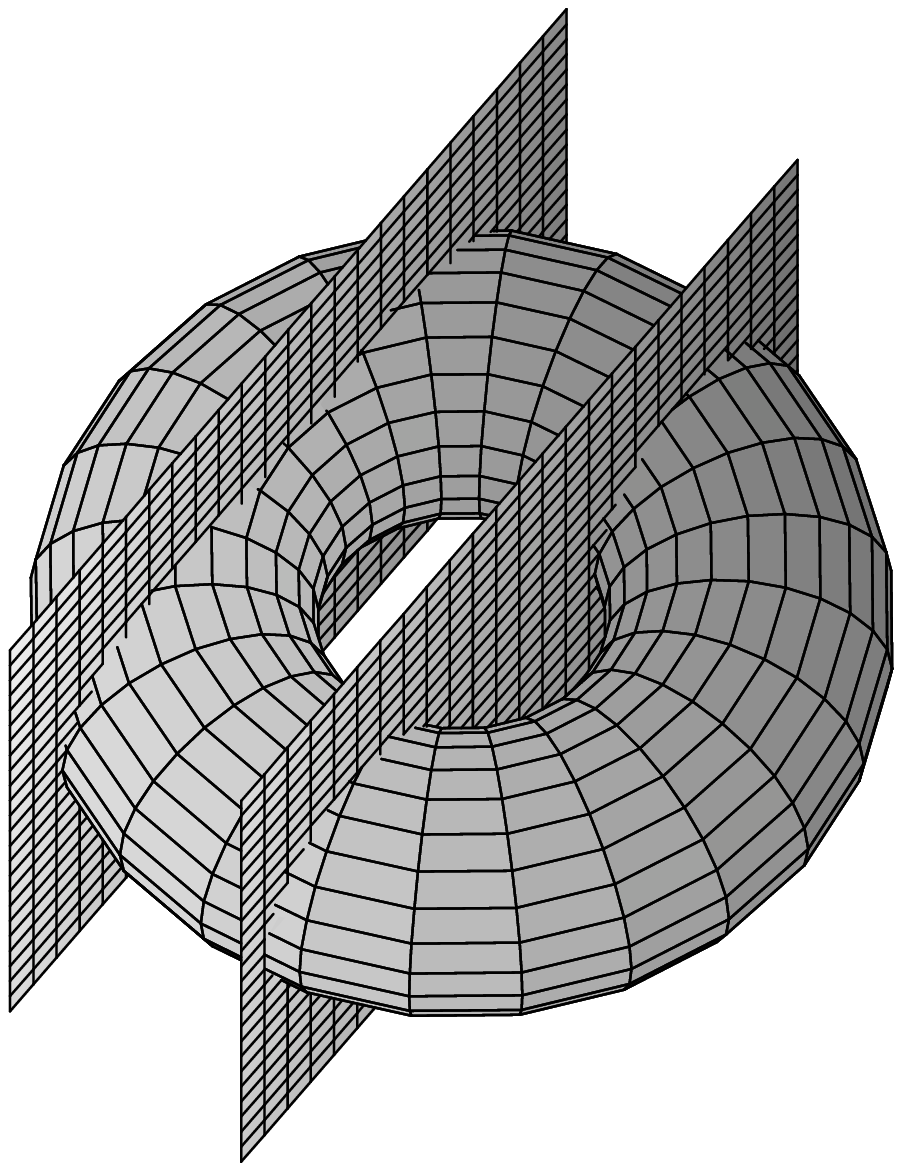}
\hspace*{0.3in}
  \epsfysize=1.2in
  \epsffile{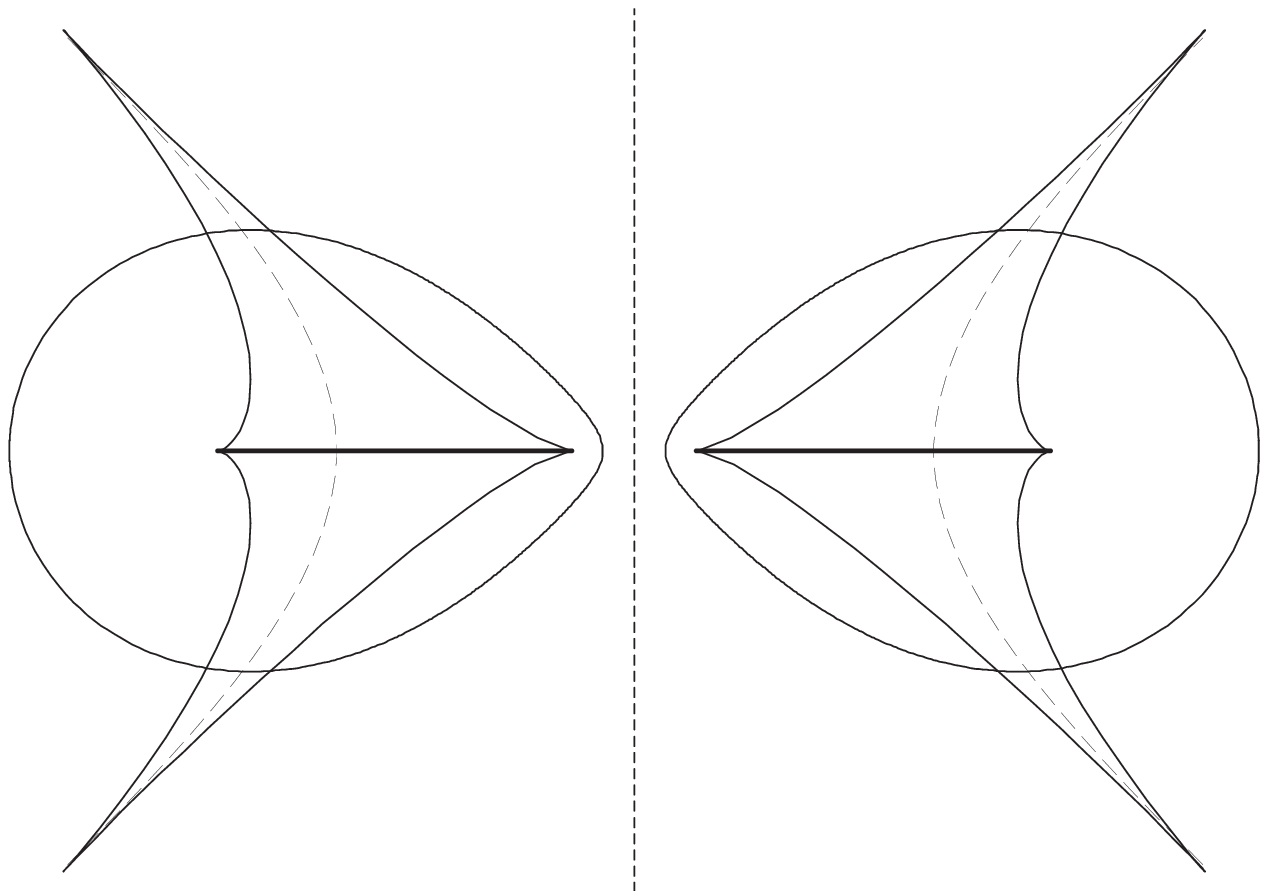}
\hspace*{0.3in}
\epsfysize=1.2in
\epsffile{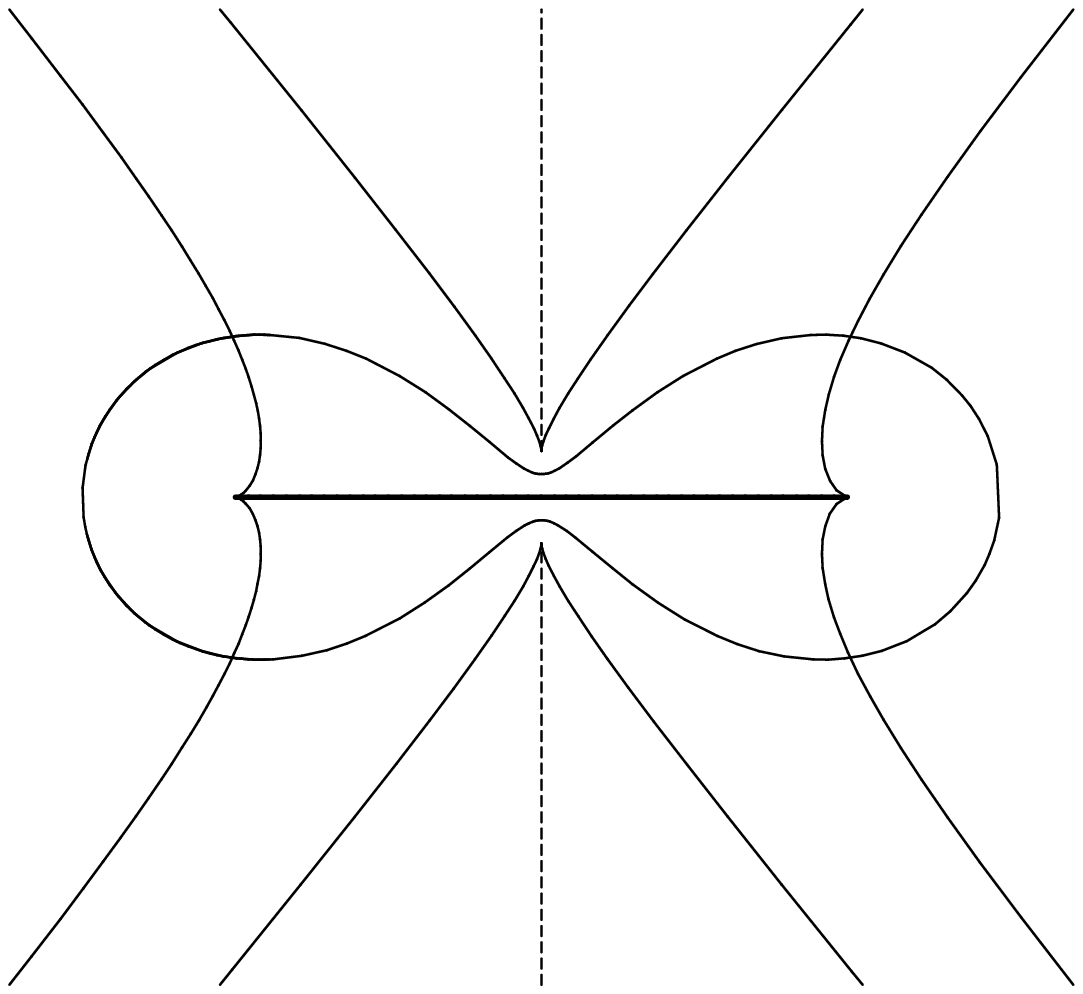}
  \end{center}
  \caption{Two plane sections of a torus close to a singular section, together with their evolutes. The thick lines are the MA and
the dashed lines are the additional parts of the SS. As the two ovals merge, two cusps on the evolute recede to
infinity, taking the branches of the SS with them. (In the right-hand figure, in fact the SS goes twice to infinity
and in between these excursions it covers the whole vertical line; this part, caused by the global structure
of the curve,  has been omitted for clarity. In this paper we
are concerned with the {\em local} behaviour of SS near the singular section.)}\label{fig:torus}
  \end{figure}

\begin{small}
{\em Acknowledgements} \
This work is a part of the DSSCV project supported by the IST Programme of the European Union (IST-2001-35443).
The authors are also grateful to Prof. Mads Nielsen and Dr. Aleksander Pukhlikov for useful discussions.
\end{small}
\section{Intrinsic geometry of generic isophote curves}\label{s:vertices-inflexions}
This section describes the geometry of isophote curves evolving on a fixed smooth surface $M$, under a 1-parameter family of
parallel plane sections. Namely, we shall examine closely the different configurations of vertices and inflexions on the sections on our surface.
  We will in particular concentrate on the evolution through a plane section which is tangent to $M$ at a point {\bf p}, so that
this section is singular. For a generic surface, three
situations arise, according to the contact between the tangent plane and $M$ at {\bf p},
as measured by the singularity type of the height function in the normal direction at {\bf p}.
See for example \cite{solid-shape} for the geometry of these situations, and \cite{b-g-t,mumford-book}
for an extensive discussion of the singularity theory.\\
$\bullet$ \ \ The contact at {\bf p} is ordinary (`$A_1$ contact'), in which case the
point is (i) elliptic or (ii) hyperbolic. The intersection of $M$ with its tangent
plane at {\bf p} is locally an isolated point or a pair of transverse arcs. \\
$\bullet$ \ \ The contact is of type $A_2$, which means that {\bf p} is parabolic. The intersection of $M$ with its tangent
 plane at {\bf p} is locally a cusped curve. \\
$\bullet$ \ \ The contact is of type $A_3$, which means that {\bf p} is a cusp of Gauss,
in which case it can be (i) an elliptic cusp, or (ii) a hyperbolic cusp. The intersection
with the tangent plane is locally an isolated point or a pair of tangential arcs.
\\
Elliptic and hyperbolic points occupy {\em regions} of $M$, separated by parabolic
{\em curves} which are generically nonsingular; on the parabolic curves are isolated {\em points} which are cusps of Gauss.

The following gives a complete description of the behaviour of vertices and inflexions  on isophotes curves near a singular point.

\begin{figure}[h!]
  \begin{center}
  \leavevmode
\epsfysize=2.2in
\epsffile{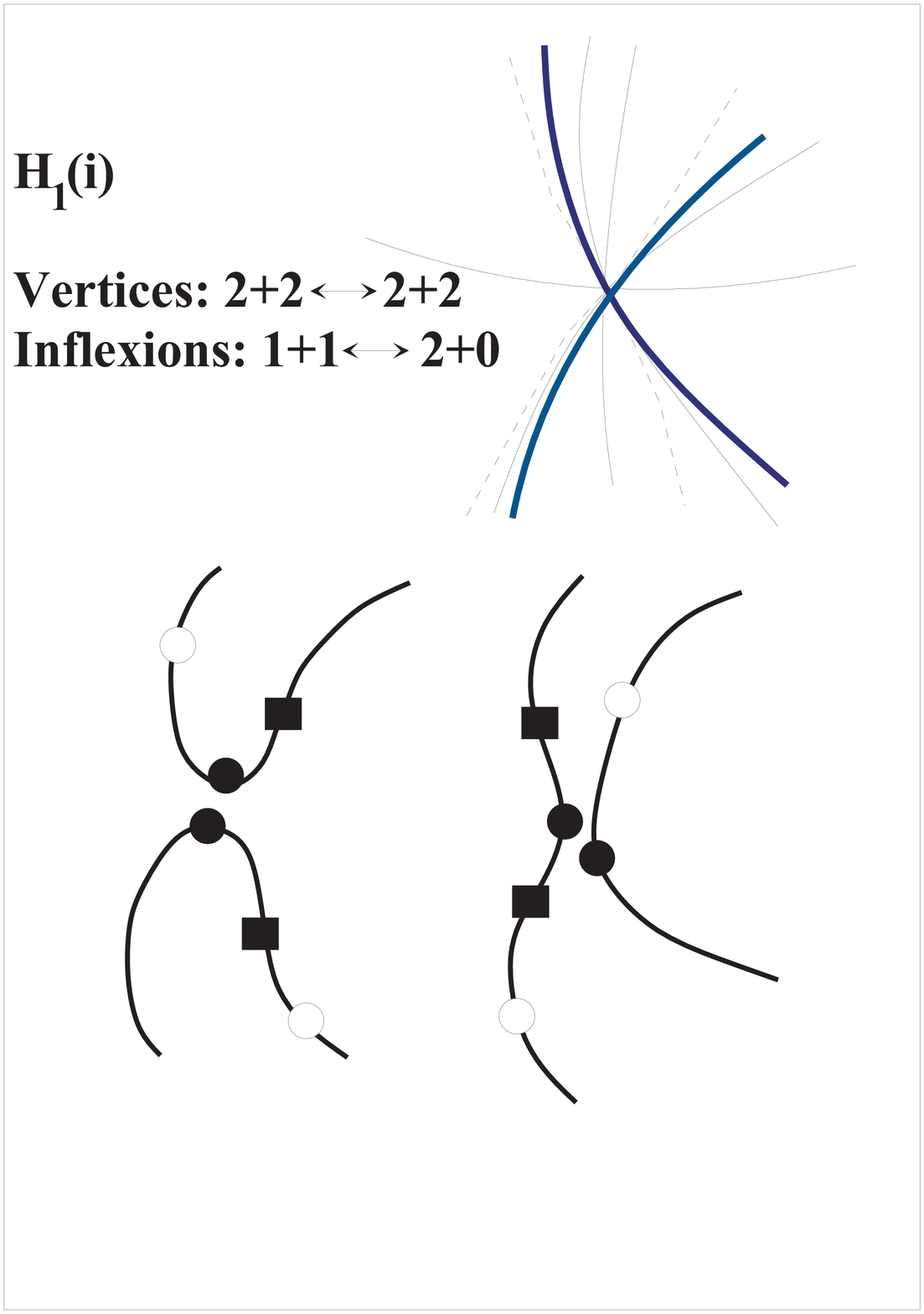}
\epsfysize=1in %
\epsffile{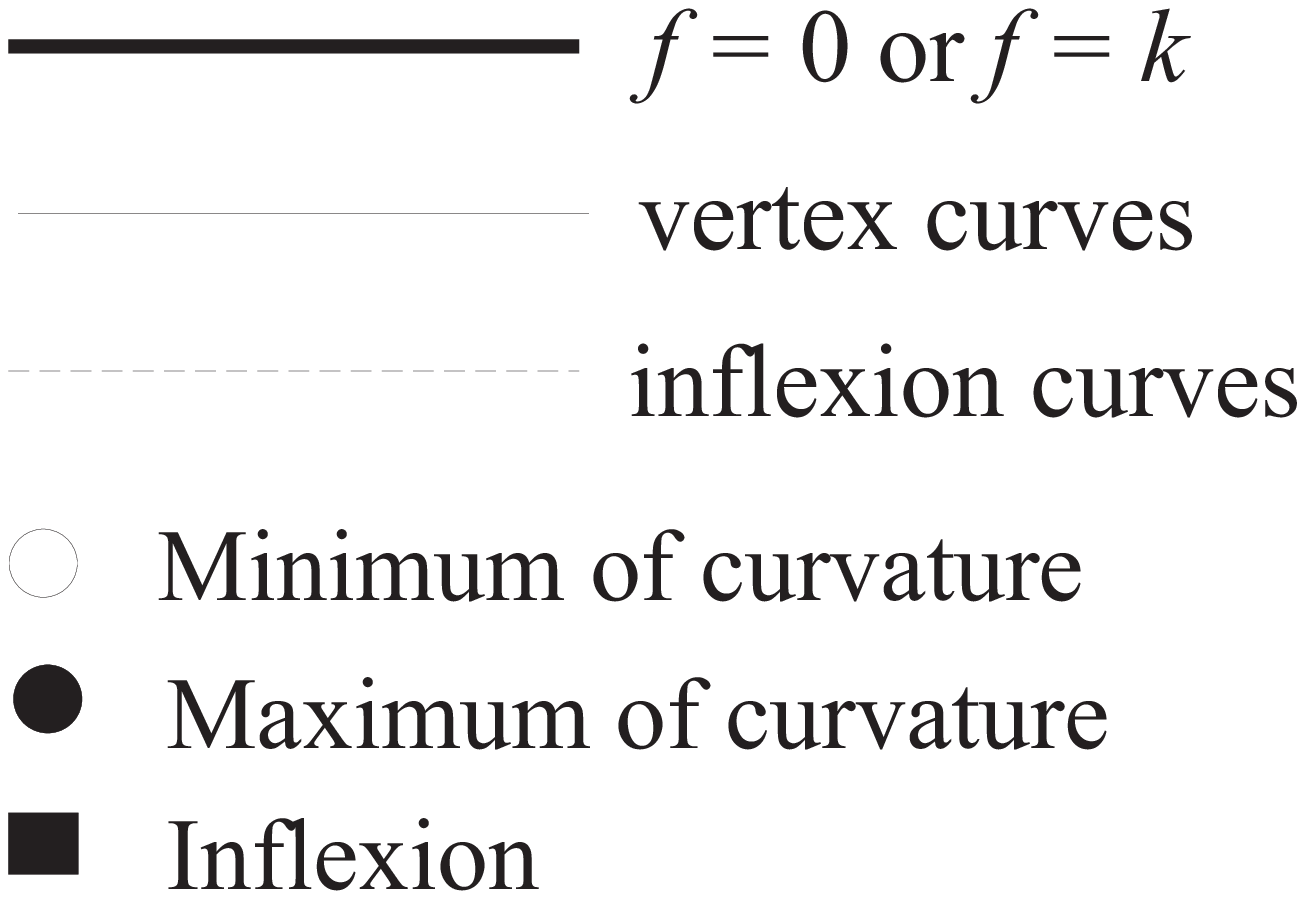} %
\epsfysize=2.2in
\epsffile{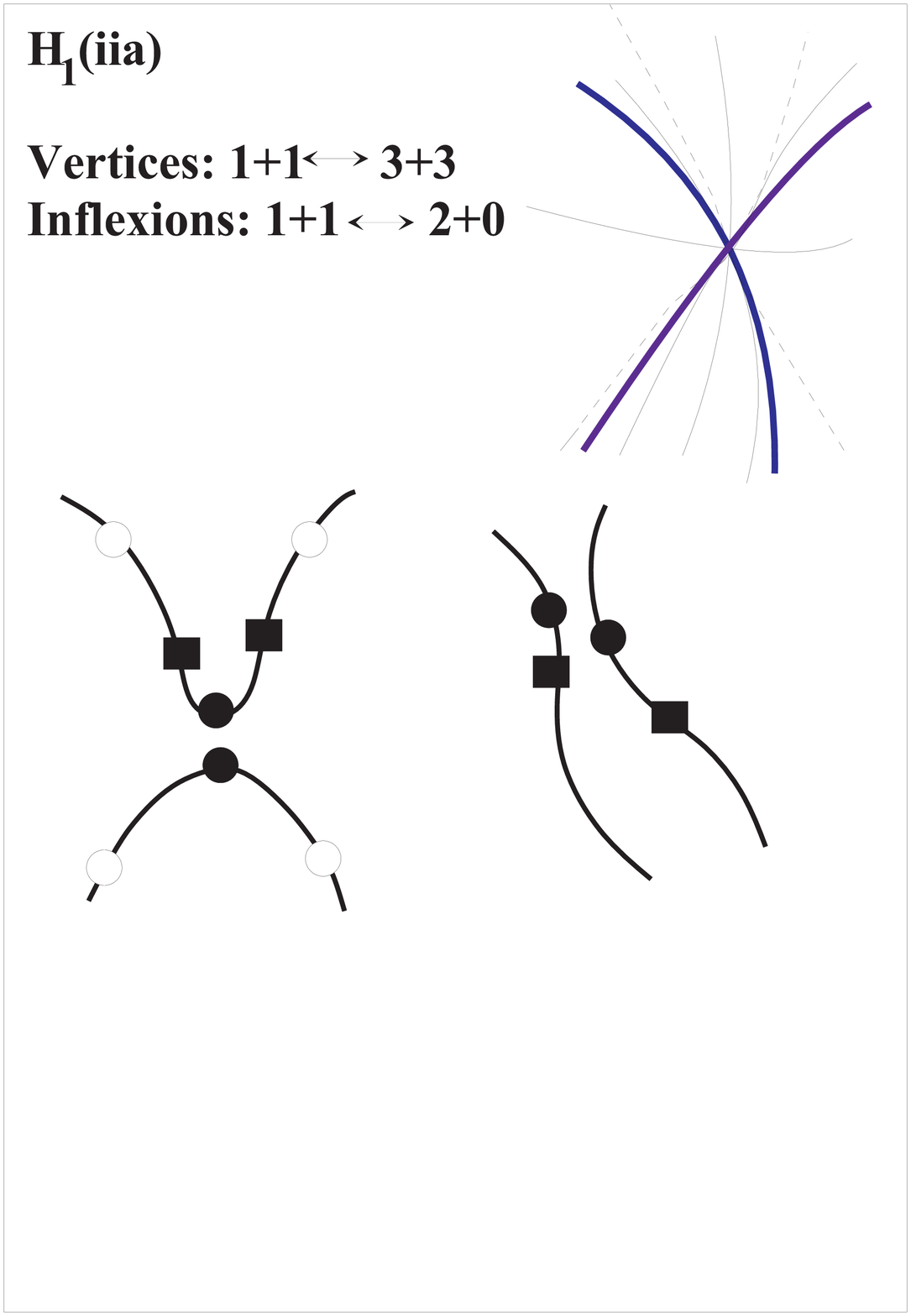}
\epsfysize=2.2in
\epsffile{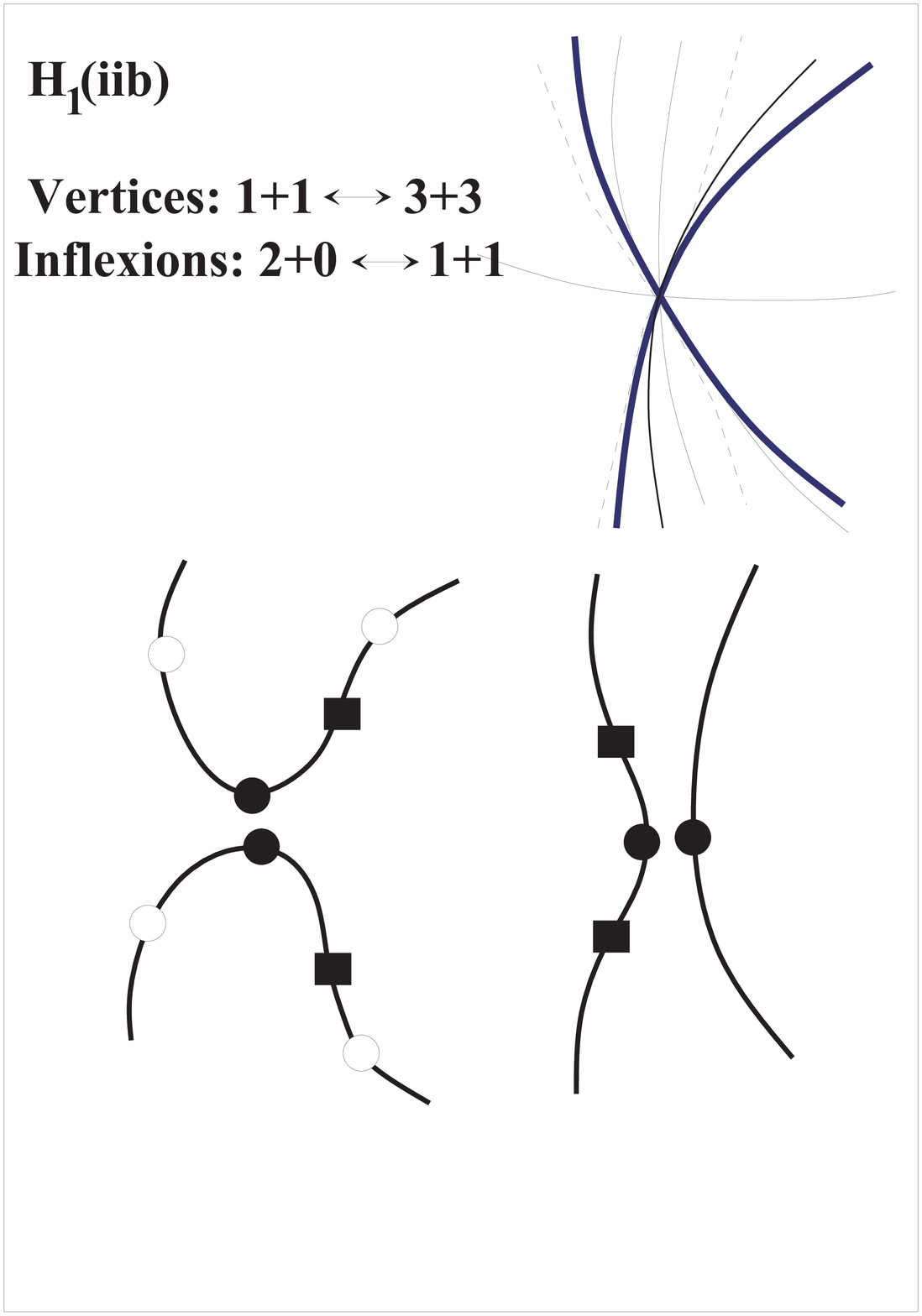}
\end{center}
  \caption{Arrangements of vertices and inflexions on the level sets of $f$, in the most generic hyperbolic case (called
${\bf H}_1$ in \cite{diatta-giblin2004}).
See Theorem~\ref{vertices-inflexions}.
In each case, we show, above, the vertex and inflexion curves---that is, the loci of vertices and inflexions on the
 level sets of $f$---and,
below, a sketch of the level curves for $f<0$, $f>0$, showing the positions of these vertices and inflexions. Thick lines:
$f=0$ or $f=k$; thin solid lines: vertex curves; dashed lines: inflexion curves. Open circles: minima of curvature; solid circles:
maxima of curvature; squares: inflexions.}
  \label{fig:hyp-geom-vertex-inflex}
  \end{figure}

\begin{theorem} \label{vertices-inflexions}
Let $f=k$ be a section of a generic surface M by a plane close to the tangent plane at  {\bf p}, $k=0$ corresponding
with the tangent plane itself. Then for every sufficiently small
open neighbourhood  $U$ of {\bf p} in M,  there exists $\varepsilon>0$ such that $f=k$ has exactly $v(\pp)$ vertices and $i(\pp)$
inflexions lying in $U$, for every $0<|k|\leq\varepsilon$, where $v(\pp)$ and $i(\pp)$ satisfy the following equalities.
\begin{enumerate}
\item[{\rm (E)}] If {\bf p} is an elliptic point, then for one sign of $k$ the section is locally empty; in the non-umbilic case,
for the sign of $k$ yielding a locally nonempty intersection we have $v(\pp)=4$, $I(\pp)=0$.
Likewise if {\bf p} is an umbilic point, then $v(\pp)=6$,  $I(\pp)=0$.
\item[{\rm (H)}]
 If {\bf p} is a hyperbolic point, v(\pp) satisfies
one of the following.
We use $\leftrightarrow$ to indicate the transition in either direction, $m+n$ indicating the numbers of vertices on the two branches of $f=k$ for
one sign of $k$ before the $\lr$ and for the other sign of $k$ after it.
In the most generic case (open regions of our surface) we have   $2+2 \leftrightarrow 2+2$ or  $1+1 \leftrightarrow 3+3$.
See Figure~\ref{fig:hyp-geom-vertex-inflex}.
In other cases, occurring along curves or at isolated points of our surface, we can have in addition
$3+2\lr 2+1$ or $3+1\lr 2+2$.
Also using the same notation, $i(\pp)$ satisfies:  $1+1 \leftrightarrow 0+2$ or  $1+2 \leftrightarrow 0+1$.
There are 8 cases in all, and the full list is given in \cite{diatta-giblin2004}.

\item[{\rm (P)}]
If {\bf p} is a parabolic point but not a cusp of Gauss, $v(\pp)=3$,  $I(\pp)=2$.
\item[{\rm (ECG)}] If {\bf p} is an elliptic cusp of Gauss, $v(\pp)=4$ ,  $I(\pp)=2$ for one sign of  $k$, and $v(\pp)=I(\pp)=0$ for the other.
\item[{\rm (HCG)}]
If {\bf p} is a hyperbolic cusp of Gauss, $v(\pp)$ satisfies $1+3\leftrightarrow 4+4$ or $2+2\leftrightarrow 4+4$,
whereas  $I(\pp)$ satisfies $2+2\leftrightarrow 0+2$, $1+1\leftrightarrow 0+0$
or $1+1\leftrightarrow 0+4$.
\end{enumerate}
\end{theorem}

For the proof and  more details see \cite{deliverables4-11}, \cite{diatta-giblin2004}.

\section{Symmetry sets (SS) and medial axes (MA) of isophote curves}
\begin{figure}[ht]
  \begin{center}
  \leavevmode
\epsfxsize=3.3in
\epsffile{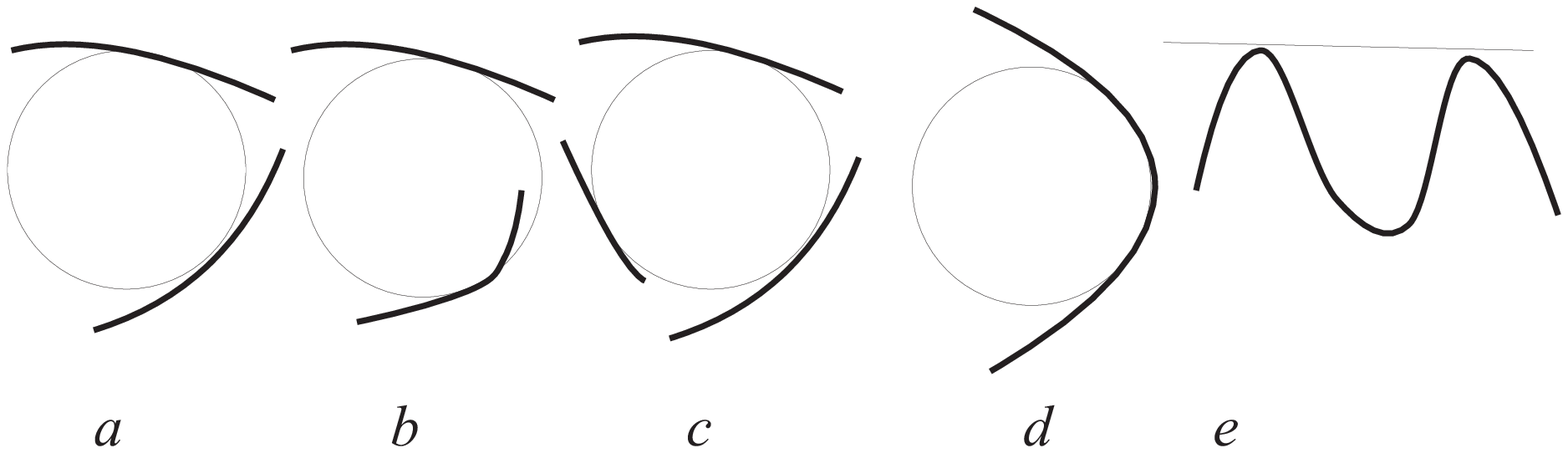}
\epsfysize=0.9in
\epsffile{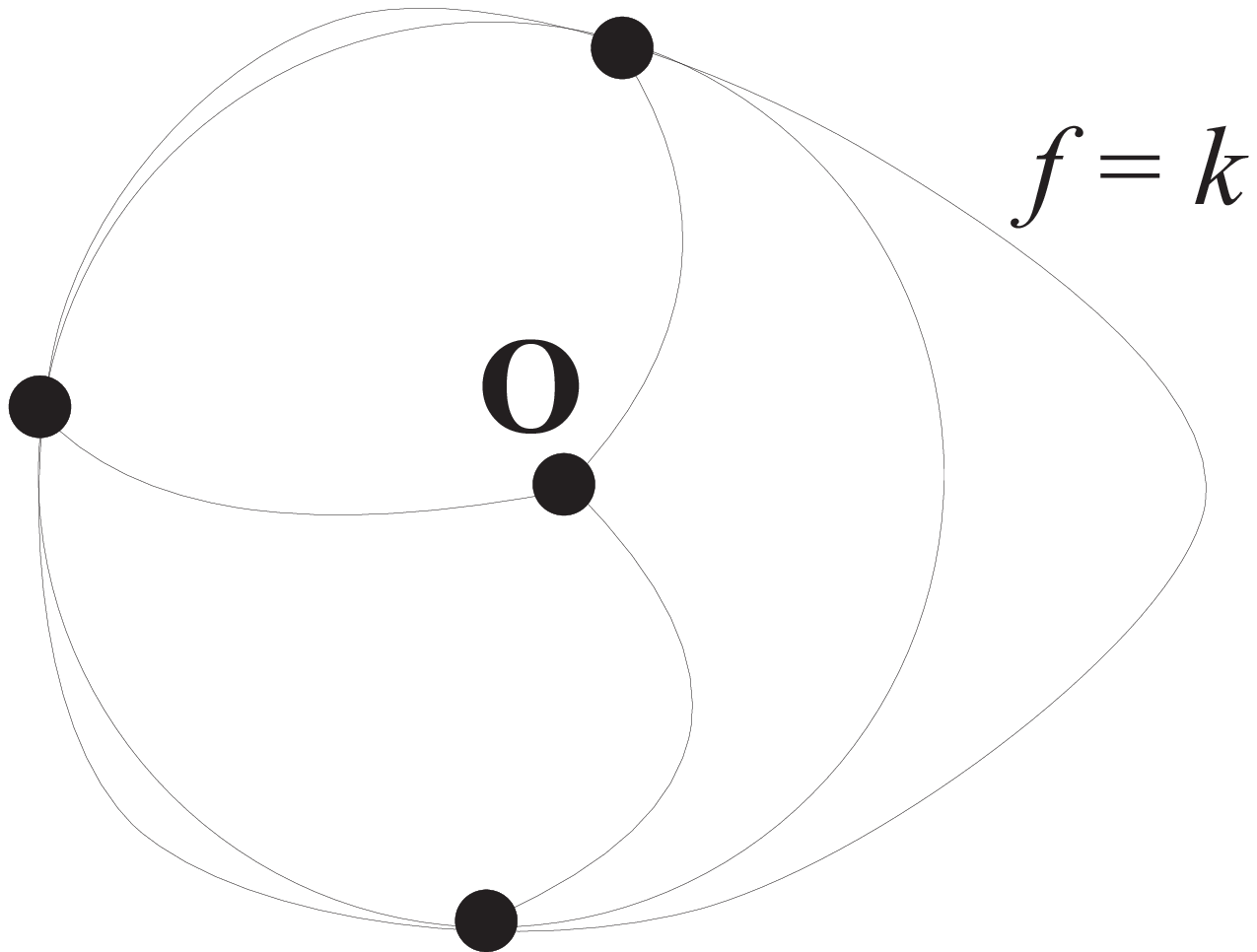}
\end{center}
  \caption{(a)-(e): Illustration of the circles whose centres contribute to the symmetry set.
(a) is an $A_1^2$, (b) an $A_1A_2$, (c) an $A_1^3$, (d) an $A_3$ and (e) a centre at $\infty$ (bitangent line).
In the last case the circle has become a straight line and the centre is at infinity. Right: a schematic drawing of a tritangent
circle and a level set $f=k$ for an umbilic point at the origin {\bf O}.
As $k\to 0$ the points of tangency trace out three curves which we call the `$A_1^3$ curves'. Calculation of these
curves is given in \S\ref{s:A13}. Once these curves are known we can calculate the locus of centres
of the tritangent circles.}
  \label{fig:order-of-contact}
  \end{figure}

The SS of a smooth simple closed curve in $\R^2$ is made of piecewise smooth curves (locus of $A_1^2$'s),
triple crossings ($A_1^3$),
cusps ($A_1A_2$), endpoints ($ A_3$) and the points at `infinity'   (they correspond to bitangent {\em lines} to the curve).
See Fig \ref{fig:order-of-contact}.

$\bullet$ $A_1^2$: The centres of bitangent circles with ordinary tangency at both points.

$\bullet$ $A_1^3$: The centres of tritangent circles with ordinary tangency at all points. They are  the triple
crossings on the symmetry set.

$\bullet$ $ A_1A_2$:  They are the centres of bitangent circles which are osculating circles at one point of the
 curve and have an ordinary
 tangency at the other point. They lie on the evolute and are cusps on the symmetry set.

$\bullet$ $ A_3$: They are the centres of circles of curvature at extrema of curvature on the curve, the endpoints
of  the symmetry set and the cusps on the evolute.

$\bullet$ Bitangent lines: the circle now has its centre at infinity so the SS goes to
infinity.

At inflexions the evolute goes to infinity and the sign of curvature changes. Thus a {\em positive maximum}
of curvature will be followed by a {\em negative minimum}, which in terms of the absolute value of curvature is again a maximum.

Our approach to the study of SS of families of curves which include a singular curve is to trace the $A_3$ points,
the inflexions, the $A_1A_2$ points and the $A_1^3$ points on the curves as they approach the point at which the
singularity develops. In this way
we obtain significant information about the SS themselves.
The patterns of vertices and inflexions have been studied in detail and for all the
relevant cases in \cite{deliverables4-11} and in \cite{diatta-giblin2004}, as recalled in Section \ref{s:vertices-inflexions}.
Subsection  \ref{s:A13} and \ref{s:A1A2} are devoted to the study of the locus of $A_1^3$ and $A_1A_2$ points,
respectively. In Subsection \ref{SS/MA} we derive information on the changes on the SS of families of isophotes curves.

\subsection{$A_1^3$ points}\label{s:A13}
The $A_1^3$ points are the centres of circles which are tangent (ordinary tangency) to $f=k$ (for any choice of
$f$, such as hyperbolic or umbilic) at three distinct points. They
occur at triple crossings on the SS. Instead of looking directly for the centres of those tritangent circles,
we rather first look for the points where those circles are tangent to the curve $f=k$ (see Fig. \ref{fig:order-of-contact}, right,
for a schematic picture of the umbilic case).   Thus we expect to have three curves, the `$A_1^3$ curves',
having the origin as their limit point, along which
the three contact points move. First, we want to find the limiting directions of these curves, ie the lines they are
 tangent to as $k\to 0$.  After finding the limiting directions, we can then
determine enough of a series expansion (possibly a Puiseux series) to decide how the $A_1^3$ curves lie with
respect to the vertex curves, etc.\ which we have determined before.  We will give an example of such a parametrization below.

The equations
which determine the $A_1^3$ curves are of course highly non-linear. They are in fact 8 equations in 9 unknowns, thereby
determining an algebraic variety in $\R^9$ which, when projected onto suitable pairs of coordinates, gives each $A_1^3$
curve in turn. There are two important features of these equations:\\
$\bullet$ Naturally they are symmetric in that the contact points can be permuted;\\
$\bullet$ The equations inevitably admit solutions obtained by making two of the tangency points coincide (`diagonal'
solutions). This causes the algebraic variety in $\R^9$ to have components of dimension greater than 1 which
we want in some way to discard.

We now set up the equations.
Any circle has the form $C(x,y)=0$ where \[C=x^2+y^2+2ax+2by+c,\] so that the centre is $(-a,-b)$ and the
radius is $r$ where $r^2=a^2+b^2-c$. However we prefer the parametrization by $(a,b,c)$ rather than $(a,b,r)$ since it results
in equations which are linear in the parameters.

Let this circle be tangent to $f=k$ at the three points $\pp_i=(x_i,y_i)$, $i=1,2,3$.  There are 8 equations $F_j=0, \ j=1,\ldots,8$ which
connect the 9 unknowns $x_i, y_i, a, b, c$.

\medskip\noindent
 $F_1:=f(x_1,y_1)-f(x_2,y_2)$,
\newline
 $F_2:=f(x_1,y_1)-f(x_3,y_3)$,
\newline
 $F_{i+2}:=x_i^2+y_i^2+2ax_i+2by_i+c, \ \ i=1,2,3$,
\newline
$F_{i+5}:= a \frac{ \partial f}{\partial y} (x_i,y_i)-b \frac{ \partial f}{\partial x} (x_i,y_i) +
x_i\frac{ \partial f}{\partial y} (x_i,y_i)-y_i\frac{ \partial f}{\partial x} (x_i,y_i), \ \ i=1,2,3.$

\medskip

The meaning of the 8 equations is as follows.

\medskip\noindent
\begin{tabular}{lll}
 $eq_1$: & $F_1=0$  &  $\pp_1$ and $\pp_2$ in the same level curve of $f$;\\
$eq_2:$ & $F_2=0$  & $\pp_1$ and $\pp_3$ in the same level curve of $f$;\\
 $eq_{i+2}:$  & $F_{i+2}=0$ &  $\pp_i$ lies on the circle $C, \ \ i=1,2,3$;\\
$eq_{i+5}$:  & $F_{i+5} =0$ & $C$ and the level set of $f$ through $\pp_i$ are tangent at $\pp_i.$
\end{tabular}

\smallskip

First from the three equations  $eq_i, \ i=3,4,5$, we can get
$a$, $b$ and $c$  as functions of $x_i, y_i$. Of course this is merely finding the circle
through three given points, which need to be non-collinear, and in particular distinct, for
a unique solution.  More details about this will appear elsewhere.

\smallskip\noindent
{\bf Remark}
In the umbilic case, we can always rotate the coordinates to make $b_0=b_2$ in the expression of $f(x,y)$, as
shown in \cite{maths-of-surfaces}. Thus, from now on we assume $b_0=b_2$ for an umbilic point.
 Once having assumed $b_0=b_2$,  we now make the genericity assumption  that $b_1\ne b_3$.
We shall also look
for solutions for these equations for which
{\em the limiting directions} (limiting angles to the positive $x$-axis) {\em are distinct}. This relates to the point made earlier, that our
equations inevitably admit `diagonal' solutions which we want to suppress.
Thus we are assuming here that the limiting directions of the three $A_1$ contact points of our tritangent circle
are distinct as the oval $f(x,y)=k$ shrinks to a point with $k\to 0$.

\begin{proposition}\label{prop:A13} Generically, there are no triple crossings, nor cusps on the local branches of the
 symmetry set of isophotes curves near a hyberbolic point.

 The limiting directions of the $A_1^3$ curves at an umbilic, making the
assumptions in the above Remark, are at angles $t_1, t_2, t_3$ equal, in some order, to
$90^\circ$, $-30^\circ$, $-150^\circ$  to the positive $x$-axis, or the `opposites' of
these, namely $-90^\circ$, $150^\circ$, $30^\circ$. This suggests strongly that there are always
two triples of $A_1^3$ contact points tending to the origin as $k\to 0$.
\label{prop:anglesA13}
\end{proposition}

Proposition \ref{prop:A13} implies, as confirmed by experimental evidence (see Fig. \ref{fig:umbilic-SS}), that there
 are in fact two triple crossings ($A_1^3$)  in the symmetry set in the umbilic case.
The proof of the Proposition is an explicit computation\footnote{This computation, like all those
underlying this article, was performed in Maple.} of the tangent cone of the algebraic variety
defined by the above equations $F_i=0, \ i=1\ldots 8$. The branches $(x_i,y_i)$
corresponding to $(t_1,t_2,t_3)=(90^\circ, -30^\circ,- 150^\circ)$ have the
form \\
$((-2b_1b_0-6b_0b_3+3c_3+c_1)t^2/6(b_3-b_1)+\ldots, t)$,\\
$(\frac{1}{2}\sqrt{3}t+\ldots, -\frac{1}{2}t + \ldots)$ and $(-\frac{1}{2}\sqrt{3}t+\ldots, -\frac{1}{2}t+\ldots)$

The actual locus of $A_1^3$ points (triple intersections) on the symmetry set close to an umbilic
point where $b_0=b_2$ as above and $b_1\ne b_3$, is $(-a(t), -b(t))$ where

$a(t)=\frac{b_0}{2}t^2 + \frac{1}{16}(7b_0b_1+9b_0b_3-3c_1-c_3)t^3 +$  h.o.t.

$b(t) = \frac{1}{8}(b_1+3b_3)t^2 + \frac{1}{16}(b_1^2+3b_1b_3+4b_0^2+5c_4-c_2-3c_0)t^3 +$  h.o.t.

Generically this curve has an {\em ordinary cusp} at the origin.
\subsection{$A_1A_2$ points}\label{s:A1A2}

 The {$A_1A_2$ points are the centres of bitangent circles which are osculating at one point and have an
ordinary tangency at the other one; they produce cusps on the symmetry set.
As in the case of $A_1^3$ points (\S\ref{s:A13}), we look in the first instance for the points where those circles are tangent to the level sets of $f$.

We find these curves by taking the circle $C$ to have equation $x^2+y^2+ax+by+c=0$ as
in \S\ref{s:A13}.  This time after elimination of $a, b, c$ we obtain 3 equations in 4 unknowns instead of 5 equations in 6 unknowns.
Let the  circle $C$ be  tangent to the same level set $f=k$ at the two points $\pp_i=(x_i,y_i), \ i=1,2$. We proceed to
write down the corresponding conditions, defining functions $F_i$ as follows.

\medskip\noindent
$F_1:= f(x_1,y_1)-f(x_2,y_2)$,
\newline
$F_2:=2a(x_1-x_2)+2b(y_1-y_2)+x_1^2+y_1^2-x_2^2-y_2^2$,
\newline
$F_3:=a f_y (x_1,y_1)-b f_ x (x_1,y_1)+x_1f_y (x_1,y_1)-y_1f_x (x_1,y_1) $,
\newline
$F_4:=a f_y (x_2,y_2)-b f_x (x_2,y_2)+x_2 f_y (x_2,y_2)-y_2 f_x (x_2,y_2)$,
\newline
$F_5:=(a+x_2)(f_{xx}f_y^2 -2f_{xy}f_x   f_y  +
f_{yy} f_{x}^2  )-f_{x} (f_x^2   +f_y^2  )$ (derivatives at $(x_2,y_2)$.

\medskip

We have the corresponding equations and their interpretations:

\begin{tabular}{lll}
$eq_1:$ & $F_1=0$ &  $\pp_1$ and $\pp_2$ are in the same level set of $f$; \\
$eq_2:$ & $F_2=0$ &a circle  with centre $(-a,-b)$ passes through $\pp_1$ and  $\pp_2$; \\
$eq_3:$ & $F_3=0$ &this circle is tangent to the level set of $f$ through $\pp_1$;\\
$eq_4:$ & $F_4=0$ &this circle is tangent to the level set of $f$ through $\pp_2$; \\
$eq_5:$ & $F_5=0$ &this circle is the circle of curvature of the level set through $\pp_2$.
\end{tabular}
\medskip
We solve $eq_2, eq_3$ for $a$ and $b$ and substitute in $eq_4$ and $eq_5$. We summarize the results as follows.
 See  Figure~\ref{fig:A1A2umbilicdirections}. We assume as before that
the limiting angles at which the $A_1$ and $A_2$ points approach the origin are distinct.

\begin{proposition}\label{prop:A1A2}  Generically, there are no cusps on the local branches of the symmetry set of isophotes
curves near a hyperbolic point.
The limiting angles in the
umbilic case must be one of the following.\\
$A_1: -30^\circ, \ A_2: 90^\circ;$ or $A_1: 150^\circ, \ A_2: -90^\circ;$\\
$A_1: -150^\circ, \ A_2: 90^\circ; $ or $A_1: 30^\circ, \ A_2: -90^\circ;$\\
$A_1:60^\circ, \ A_2: -120^\circ$ or vice versa; \\
$A_1:-60^\circ, \ A_2: 120^\circ$ or vice versa.
\label{prop:anglesA1A2}
\end{proposition}

\begin{figure}
\begin{center}
\includegraphics[width=5in]{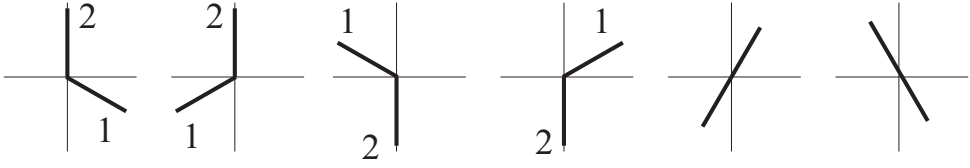}
\end{center}
\caption{The possible limiting directions of $A_1$ and $A_2$ contact points of $A_1A_2$ circles, in the umbilic
case with axes rotated so that $b_0=b_2$, assuming the limiting directions are unequal.  Those labelled 1
are only $A_1$ directions and similarly for
$A_2$. The unlabelled directions can be either, with $A_1$ and $A_2$ at $180^\circ$ to one another.}
\label{fig:A1A2umbilicdirections}
\end{figure}

This means that there are six cusps ($A_1A_2$) on the SS in this umbilic case. In that case,
we expect each cusp (which requires  an $A_1$ and an $A_2$ contact) to use one of
the above six solutions, for a definite choice of $A_1$ and $A_2$ in the last two cases.

\subsection{Symmetry Sets (SS) and Medial Axes (MA)}\label{SS/MA}
As suggested by Theorem \ref{vertices-inflexions}, Propositions \ref{prop:A1A2} and  \ref{prop:A1A2}, the local structure of the SS and MA
of individual isophote curves and its transitions are as follows:

$\bullet$ parabolic points: the local structure of SS is just $3$ separate branches correponding to the $3$ vertices separated by inflexions (Theorem \ref{vertices-inflexions}), see Fig. \ref{fig:transitionparabolic}.

$\bullet$ nonumbilic elliptic points: the SS is made of just $2$ transverse arcs, one joining two centres of curvature at maxima of curvature and the other one two minima of curvature. The SS will look like itself and disappear as the curve shrinks to a point.
\begin{figure}
  \begin{center}
 \leavevmode
\epsfysize=2in
\epsffile{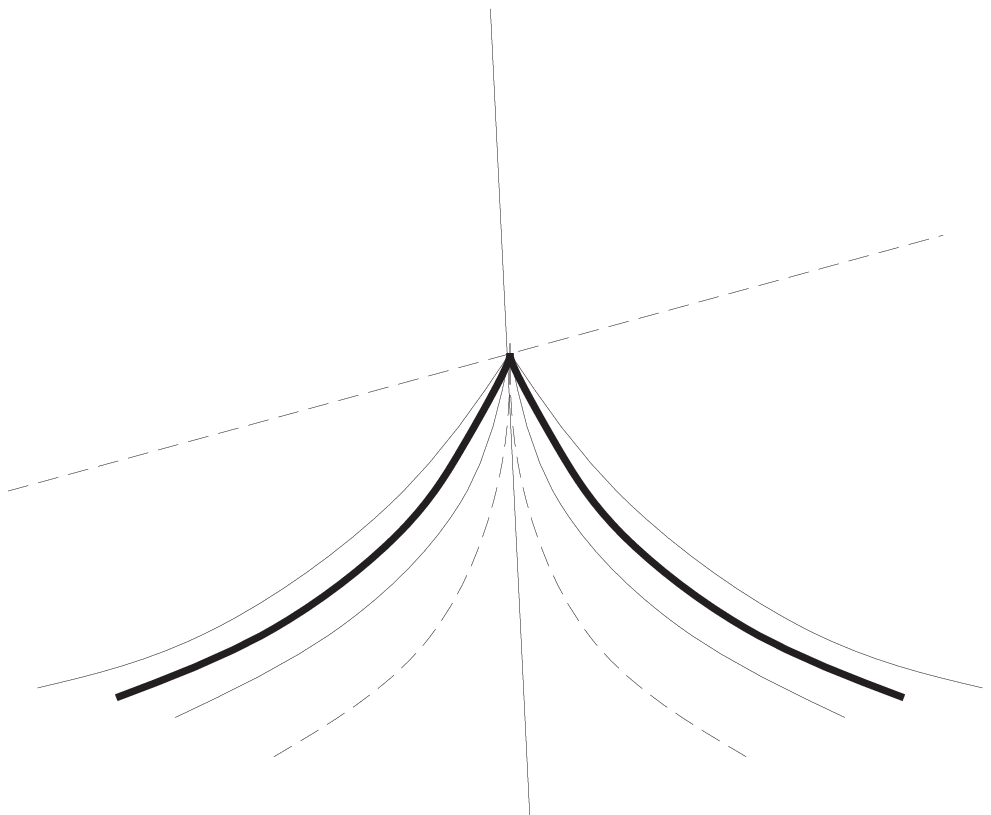}
\epsfysize=1.7in
\epsffile{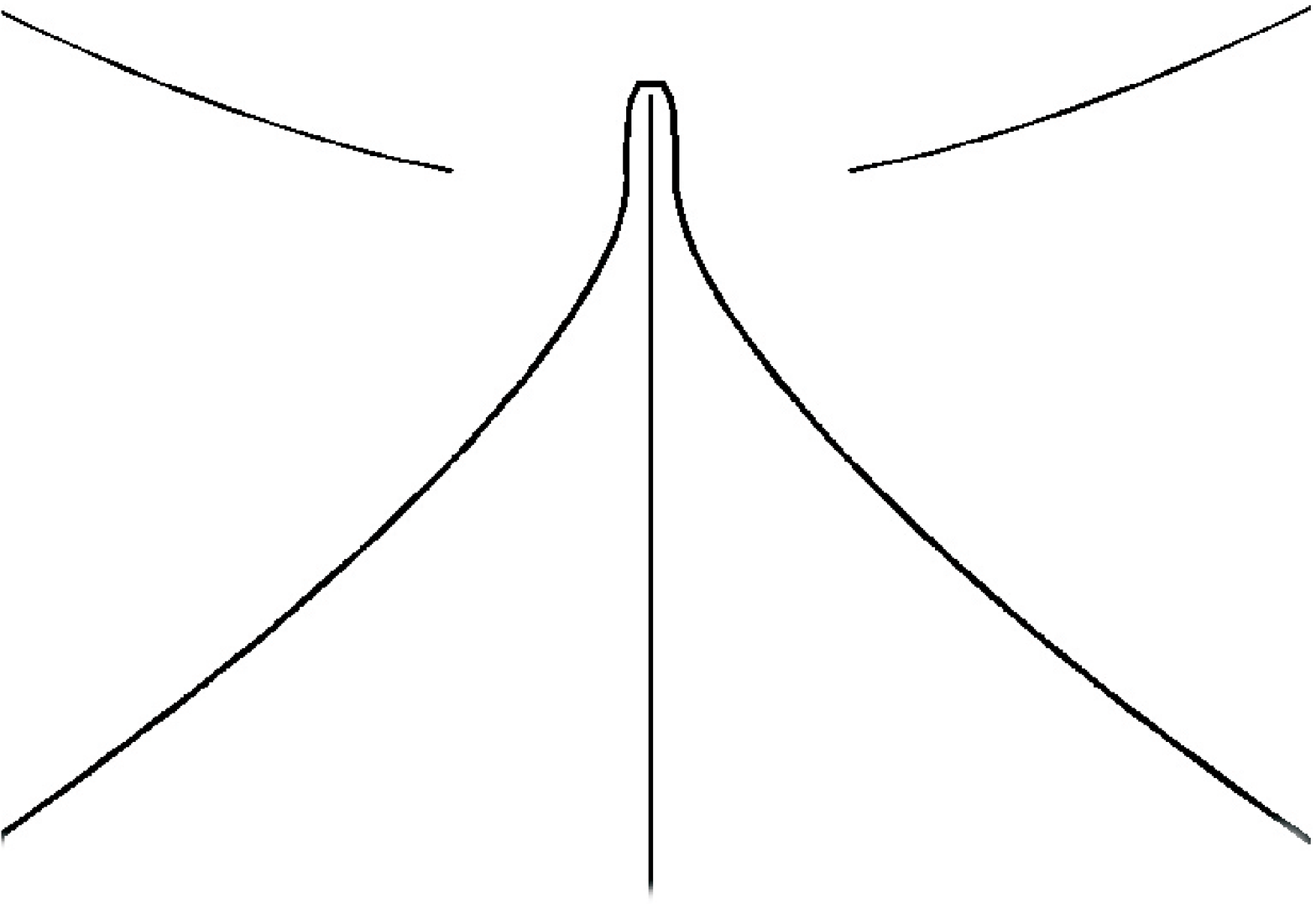}
\epsfysize=1.7in
\epsffile{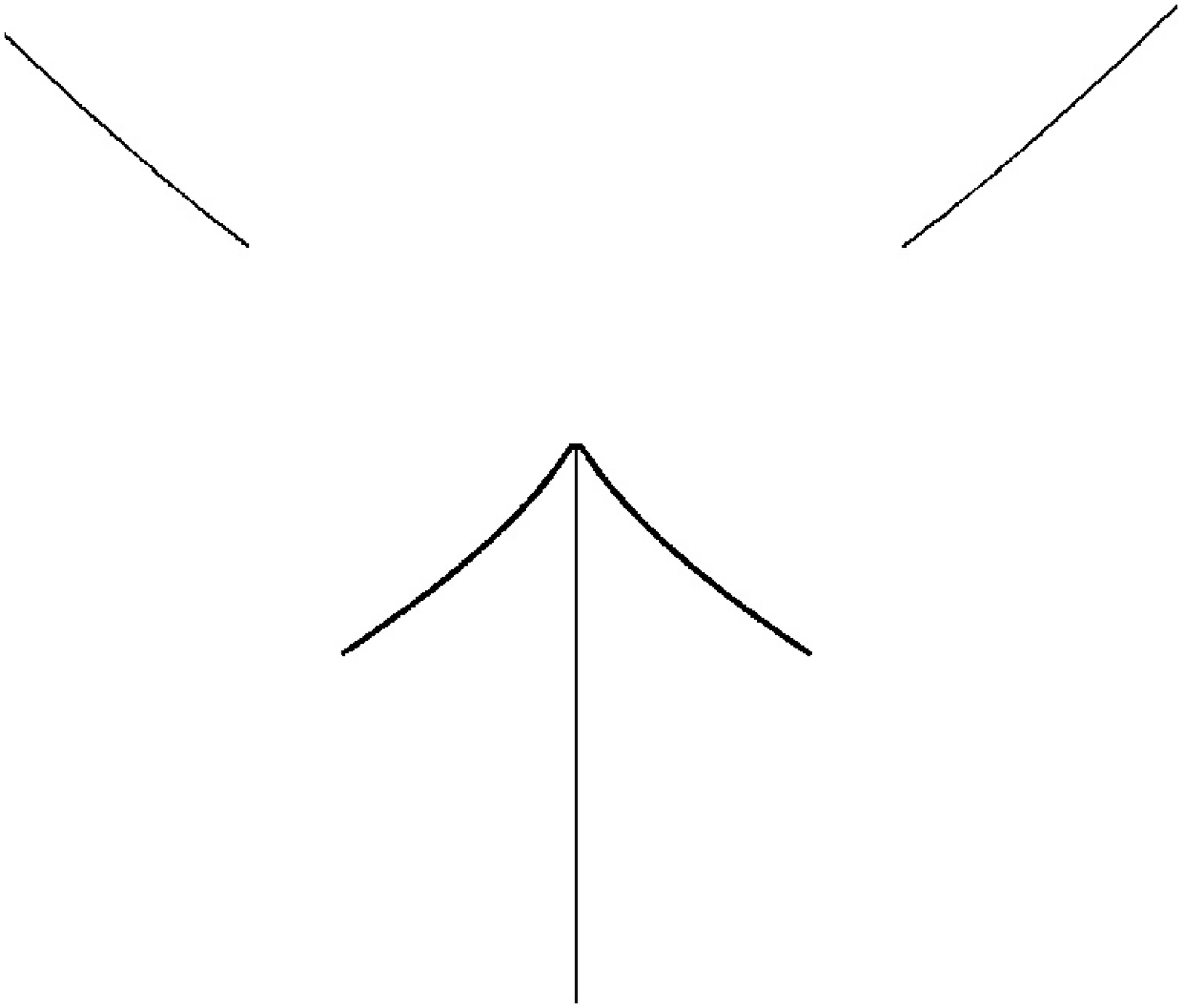}
  \end{center}
\caption{To the left:  A schematic picture of the  patterns of the vertices (vertex set $V_p=0$: thin solid line) and  inflexions
 (inflexion set $I_p=0$: dashed line) of the level curves $f=k$ evolving
through  a parabolic point, together with the zero level set $f=0$ (thick line).
 The vertex set has two cuspidal branches and one
smooth branch.  The inflexion set has one cuspidal branch which is always below all cupidal branches and one smooth branch.
The zero level set $f_p=0$ has one cuspidal branch which is always between the two cuspidal vertex branches.
The level set $f_p=k$ then evolves so that the number of vertices remains as 3 and the number of inflexions
as 2 for both signs of $k$, with $k$ small.
 In the middle and to the right: Symmetry sets (thin lines) of curves
(thick lines) which are sections of a surface close to the tangent plane at a parabolic point. One sees $3$ vertices sepated by two inflexions on both
before and after the transition. At the transitional moment itself, the branches reach
 right to the curve, which then has an ordinary cusp. Figure produced with LSMP\protect{\cite{lsmp}}.}
\label{fig:transitionparabolic}
  \end{figure}

$\bullet$ hyperbolic point: the SS and MA are made of smooth branches, which do not connect together to form cusps or crossings.
This implies in particular that generically, the SS (and MA) is just given by the geometry of vertices and inflexions
 as well as how they are distributed along the isophote curves, as described in Section \ref{s:vertices-inflexions}. The branches of the SS
 will start at endpoints which are the centres of curvature of the isophote curves at vertices and they point towards the
corresponding vertex if the isophote curve has a local minimum of curvature, and away from the vertex
where the curve has a maximum of curvature.

$\bullet$ near umbilics: the SS have generically two triple crossings and six cusps.
 Hence generically,  the SS has one structure, as in Fig. \ref{fig:umbilic-SS}.

\begin{figure}[h!]
  \begin{center}
  \leavevmode
\epsfysize=3in
\epsffile{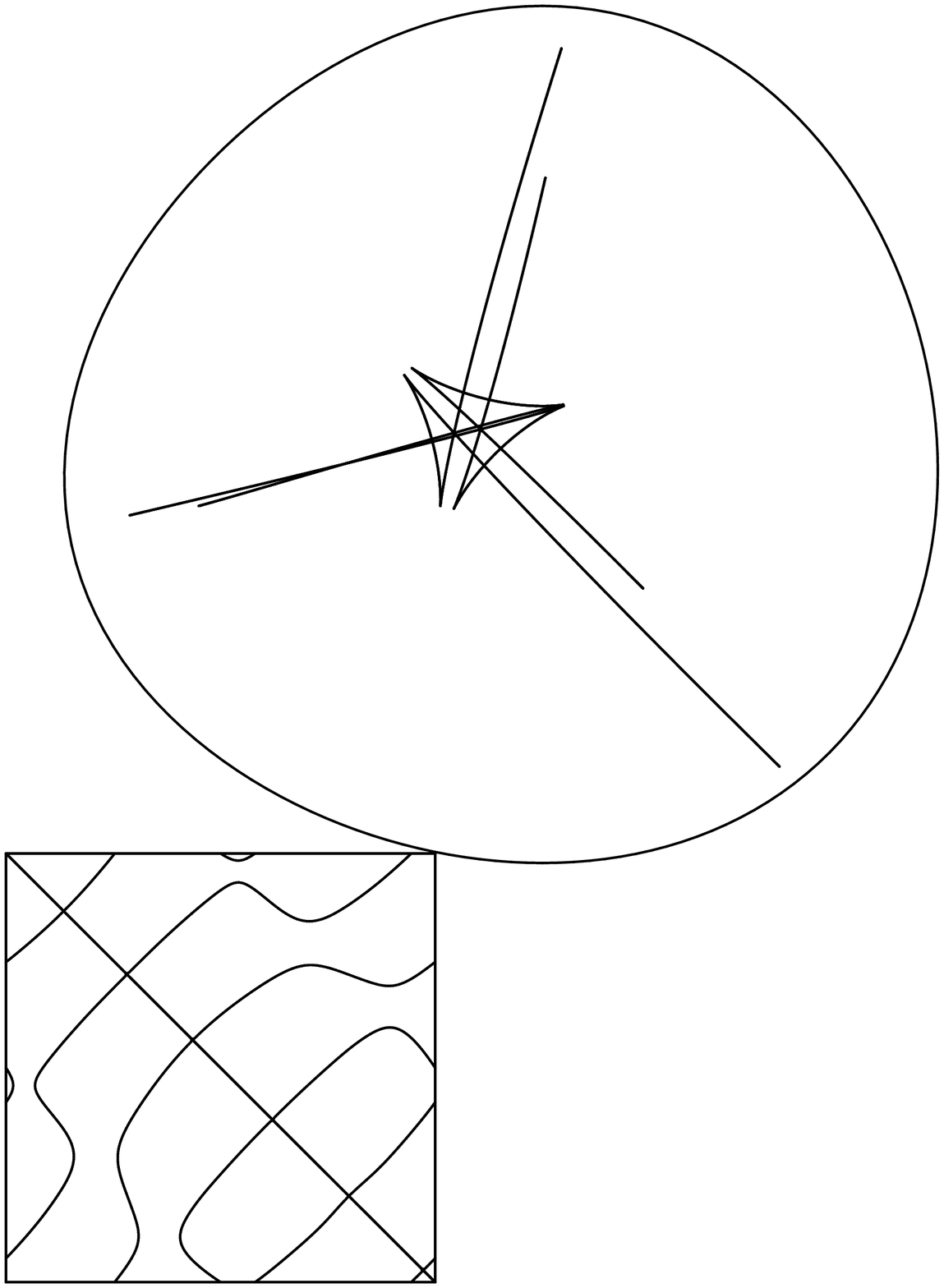}
\epsfysize=3in
\epsffile{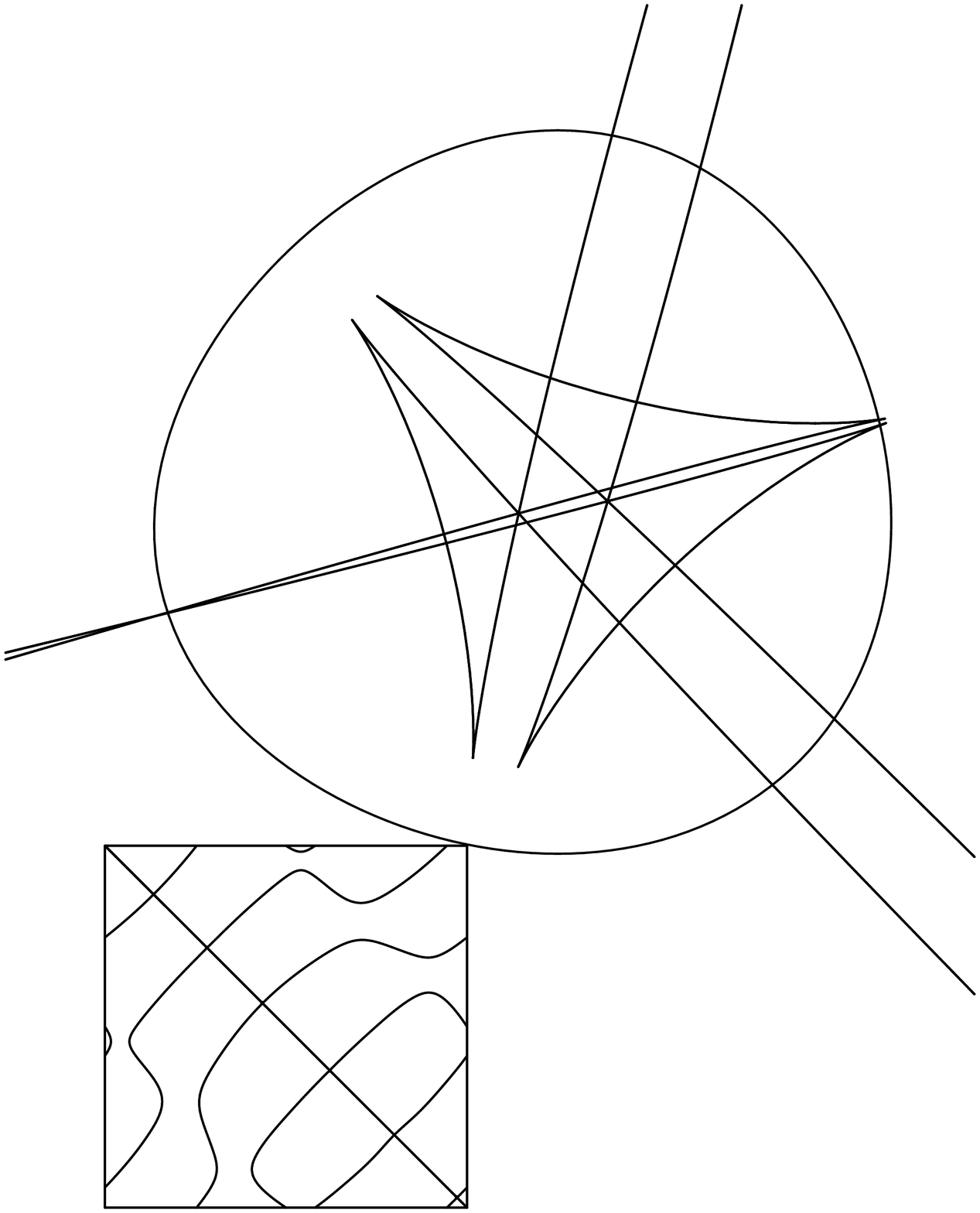}
   \end{center}
  \caption{Symmetry set and pre-symmetry set  of $f=k$,  in the umbilic case $f(x,y)= x^2+y^2+x^3-xy^2+2y^3$
and $k=0.09$. The figure to
the right is the same as the left hand side one, but the symmetry set has been enlarged so that the  two $A_1^3$ points
(triple crossings)  and the six $A_1A_2$  (cusps) are more visible.  One can also see the six endpoints of the symmetry set,
 corresponding to the six vertices on the curve. Varying $k$ then the SS will still look like itself and disappear as $k \to 0$.
This figure illustrates the results of Proposition~\ref{prop:anglesA13}
and Proposition~\ref{prop:anglesA1A2}.}
  \label{fig:umbilic-SS}
  \end{figure}
For the drawing of the SS and MA, we will need the pre-symmetry set (preSS)
which is the subset of the cartesian product $I\times I$ of the parameter space $I$, defined by  the pairs $(s,t)$
 corresponding to points ${\text \bf p}=\gamma(s)$ and ${\text \bf q}=\gamma(t)$ which contribute to the SS. That is,
there is a circle tangent to $\gamma$ at the points $\gamma(s)$ and $\gamma(t)$. See Fig. \ref{fig:umbilic-SS}

\section{Evolution of symmetry sets of isophote curves in $1$-parameter families of surfaces}
As explained in Section \ref{s:vertices-inflexions}, given a generic surface $M$, elliptic and hyperbolic points occupy
{\em regions} of $M$, separated by parabolic {\em curves} with  isolated {\em points} on them which are cusps of Gauss.
We can then consider moving from a hyperbolic point to a parabolic point of $M$.  We can also realise this
by evolving the surface in a 1-parameter family,
of the form $z=x^2-\alpha^2y+b_0x^3+b_1x^2y+b_2xy^2+b_3y^3+...$, where
$\alpha\to 0$ and $b_3\neq 0$.
 It turns out that, generically, the only hyperbolic points which exist sufficiently near a parabolic point are the ones
corresponding to vertex transition $1+1\lr 3+3$ in Theorem \ref{vertices-inflexions}.
The Figure \ref{fig:transition-hyp-para-vert} shows how the vertices  behave on a 2-parameter family of plane sections
near the  tangent plane at a hyperbolic point, evolving to a family of sections near a parabolic point.

\begin{figure}[h!]
  \begin{center}
  \leavevmode
\epsfysize=1.7in
\epsffile{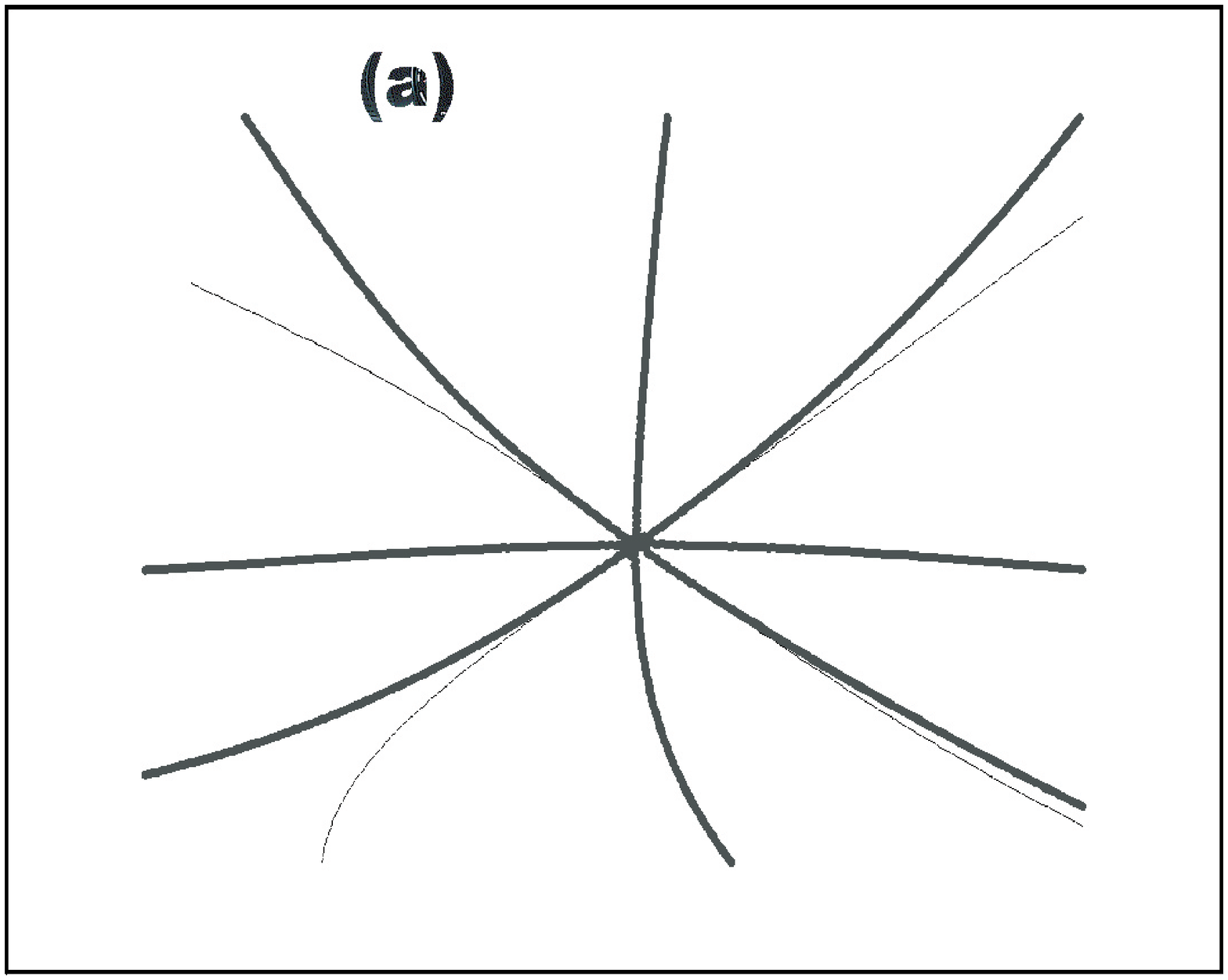}
\epsfysize=1.7in
\epsffile{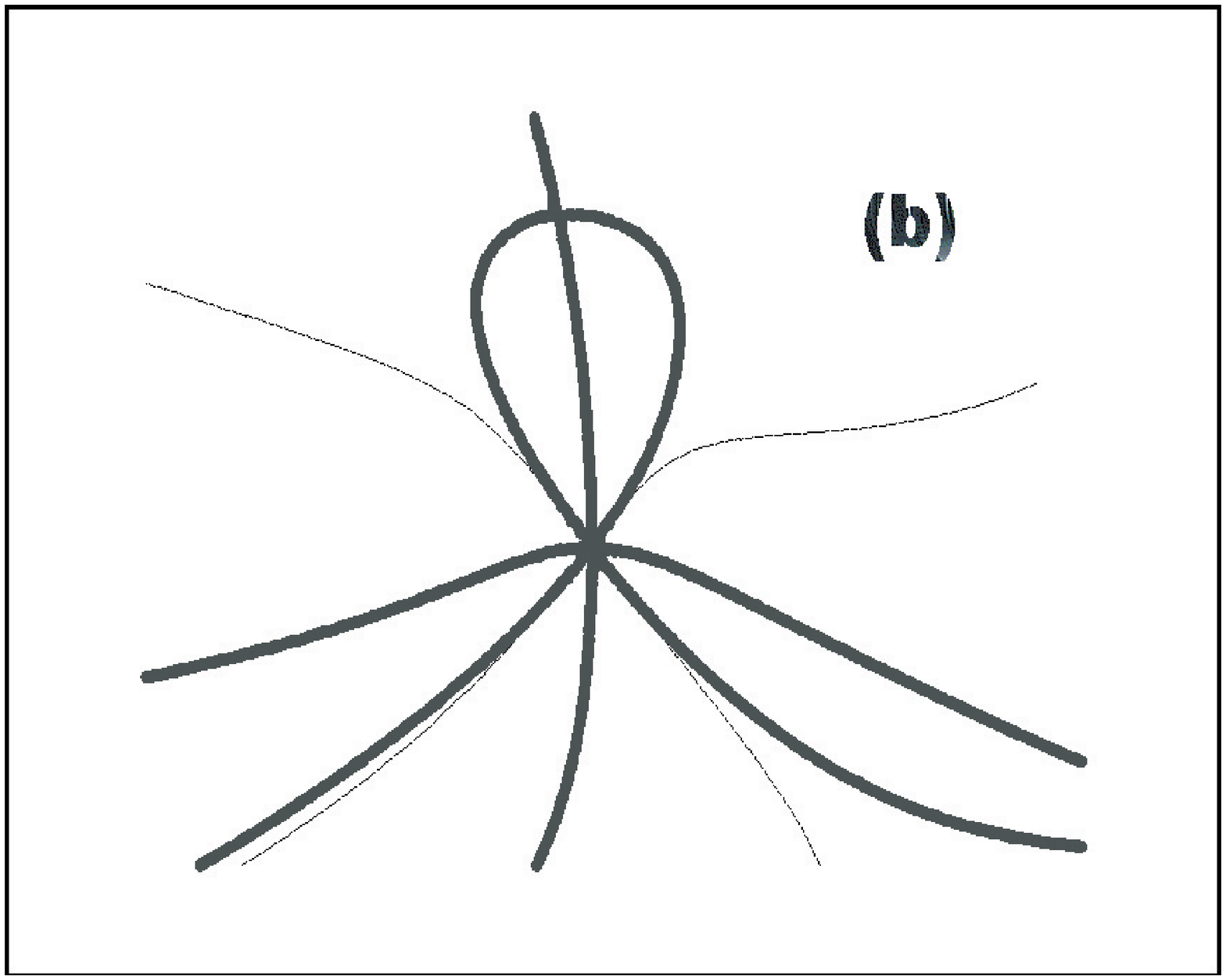}
\epsfysize=1.7in
\epsffile{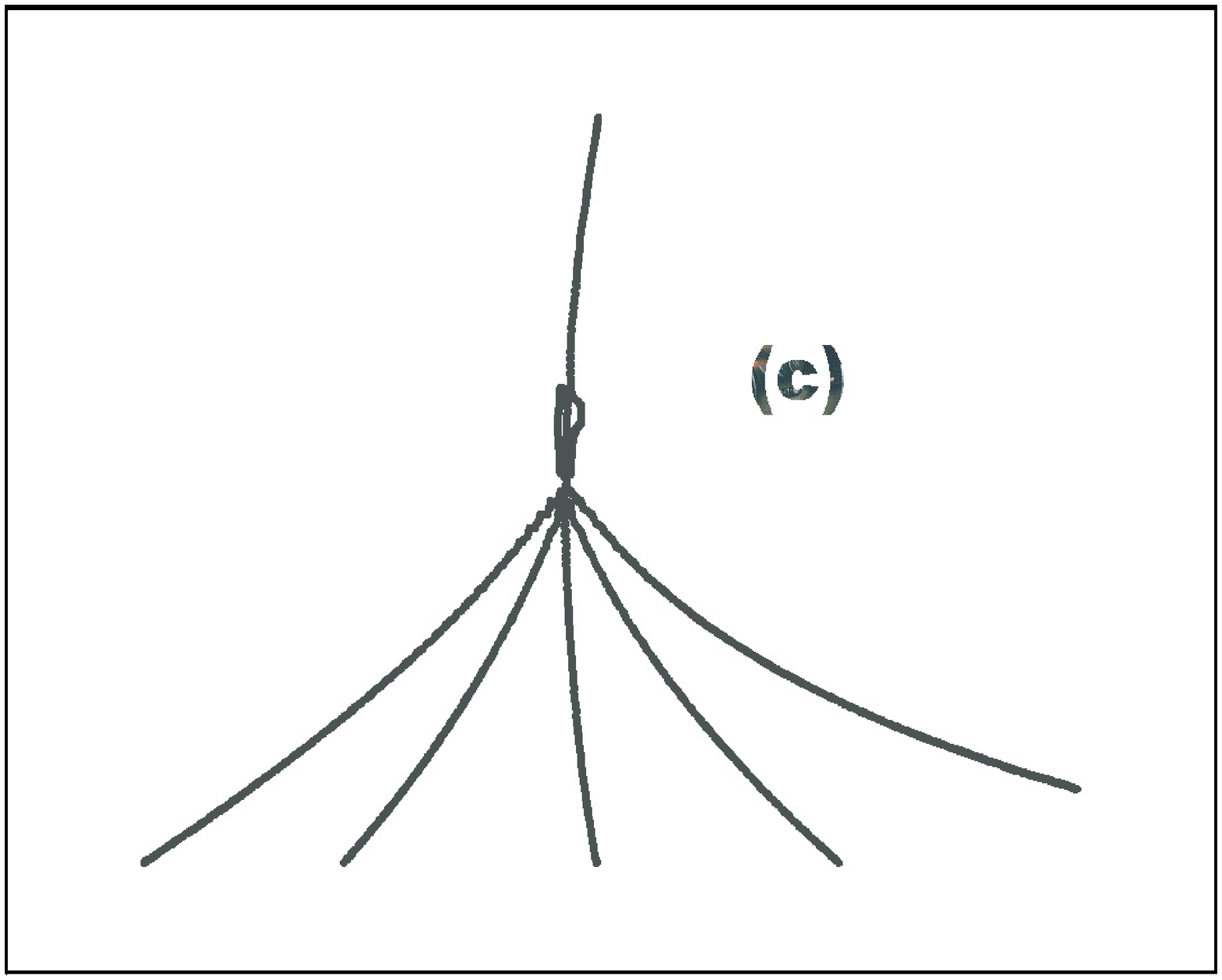}
\epsfysize=1.7in
\epsffile{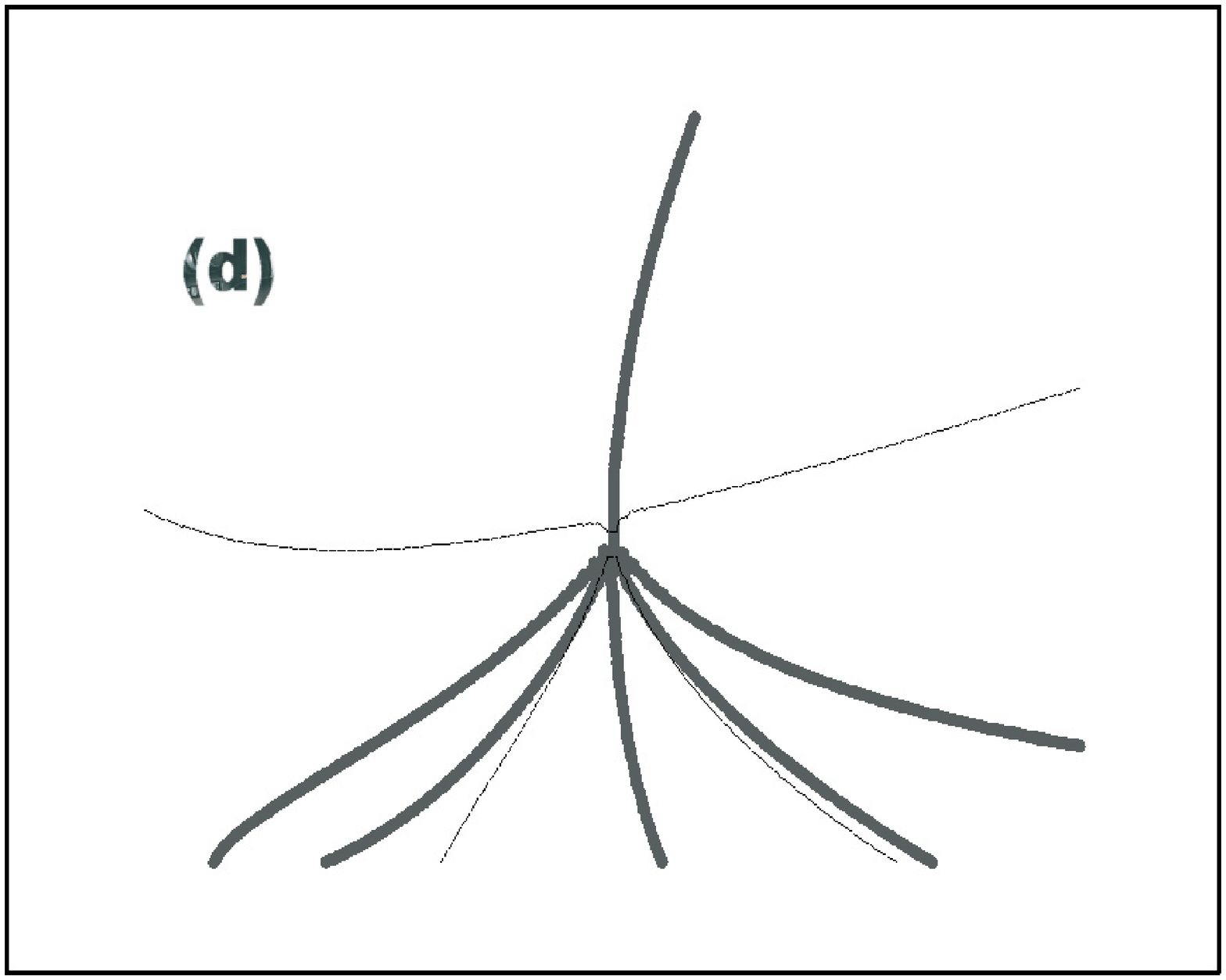}
\end{center}
  \caption{Transition of the patterns of vertices (thick curve) and inflexions (fin curve) on a 2-parameter family of plane sections $f_\alpha=k$ near the  tangent plane at a
hyperbolic point, evolving to a family of sections $f_0=k$ near a parabolic point. $f(x,y)=x^2-\alpha^2y^2+x^3+2x^2y-xy^2+y^3$. (a) $\alpha=1$ (hyperbolic);  the vertex set has $4$ branches and the inflexion $2$.
(b) $\alpha=0.3$, the top part (above x-axis) of two vertex branches join together to form a loop which is shrinking to the origine as $\alpha\to 0$. The vertex branch tangent to $x=0$ stays smooth.
The other vertex branch bend to become a cusp. (c) $\alpha=0.05$: the vanishing loop.  (d) As $\alpha\to 0$,  the inflexion set exchanges branches: the  top part (above x-axis) join together to make a smooth branch, whereas the bottom part form a cusp below the cuspidal vertex branches. Compare Figures \ref{fig:hyp-geom-vertex-inflex}
 and
\ref{fig:transitionparabolic}.}
  \label{fig:transition-hyp-para-vert}
  \end{figure}
\section{Conclusion}
This paper represents a step towards understanding the evolution of SS and MA
of families of isophote curves, or more generally of families of plane sections of s generic surface
in 3-space.
 The evolution of the MA depends, in an essential way, upon the underlying evolution of the SS
\cite{giblin2000}, which is why we have concentrated on the SS in this paper. An
interesting follow up of this work, would be to combine into a more global represention of an image by the
collection of those individual representations, as a singular surface in scale space, whose sections are the
individual SS  and MA.

\end{document}